# A Note on Near Subnormal Weighted Shifts*

**WANG Gong-bao（王公宝）**

*(Department of Basic Courses, Naval University of Engineering, Wuhan，430033)*

**MA Ji-pu (马吉溥)**

*( Department of Mathematics, Nanjing University, Nanjing,   210093)*

**Abstract: A complete characterization of near subnormality for bilateral weighted shifts** is obtained. As an application of the main results, many new answers to the Hilbert space problem 160 are presented at the end of the paper.

**MR(2000) subject classification:**    47B37, 47B20.

**CLC number:** O177.1

## 1   Introduction

Let *H* be a Hilbert space and *B(H)* be the Banach algebra of all bounded linear operators on *H*.

In [1], two operator classes of *D*-near subnormal operators and near subnormal operators were introduced. Some properties of such operators were obtained by using the M-P generalized inverse. In addition, a new necessary and sufficient condition for an operator to be subnormal was provided.

For convenience of readers, some relevant concepts and results of [1] are cited as follows.

**Definition 1** Let  $A, D \in B(H)$  and  $D \geq 0$. If there exists a constant  $m > 0$ such that $D \geq mA^*DA,$  then *A* is called a *D*-near subnormal operator.

**Definition 2** Let *A* be a hyponormal operator on a Hilbert space *H* and

* Research supported by the NNSF(10271053) of China and the Science Foundation of Naval University of Engineering (HGDJJ03001).

$Q_A = A^*A - AA^* (Q_A \geq 0)$. If $A$ is a $Q_A$ – near subnormal operator, then $A$ is called a near subnormal operator.

**Theorem**[1] Suppose $A$ is a hyponormal operator, then $A$ is near subnormal if and only if (1) $N(Q_A) \in LatA$ and (2) $Q_A^{\frac{1}{2}} A Q_A^{+\frac{1}{2}} \in B(H)$.

In the theorem above, $N(Q_A), LatA$ and $Q_A^{+\frac{1}{2}}$ denote the null space of $Q_A$, the invariant subspace lattice of $A$ and the M-P generalized inverse of $Q_A^{\frac{1}{2}}$ respectively(cf.[2]).

From [1], we know that the class of subnormal operators is a proper subset of near subnormal operators. And the class of near subnormal operators is a proper subset of hyponormal operators. In [3], Peng presented a necessary and sufficient condition for unilateral weighted shifts to be subnormal. In [4], we obtained necessary and sufficient conditions for unilateral weighted shifts to be near subnormal. Moreover, many answers to the Hilbert space problem 160 (cf.[5]) were provided .

The main purpose of this paper aims at discovering which bilateral weighted shifts are near subnormal. Necessary and sufficient conditions are obtained in terms of the weight sequence of a bilateral weighted shift. In the 160th problem of [5], it is required to construct an operator which is hyponormal but not subnormal. Halmos regarded this problem as non-trivial. As an application of the main results of the paper, we'll present many new answers to the 160th problem of [5] at the end of this paper.

## 2  Main Results and Proof

Throughout the following sections, let $H$ be a separable complex Hilbert space with an orthonormal basis $\{e_n\}_{n=-\infty}^{+\infty}$, and $T \in B(H)$. $T$ is called a bilateral weighted shift if there exists a bounded sequence of complex numbers $\{\beta_n\}_{n=-\infty}^{+\infty}$ such that $Te_n = \beta_n e_{n+1}$ for all $n = 0, \pm 1, \pm 2, \cdots$. In this case, $\{\beta_n\}_{n=-\infty}^{+\infty}$ is called the weight sequence of $T$. It is easy to see that $T$ is hyponormal if and only if $|\beta_n| \leq |\beta_{n+1}|$ for all $n = 0, \pm 1, \pm 2, \cdots$. And it is evident that $T^*e_n = \overline{\beta_{n-1}} e_{n-1}$ for all $n = 0, \pm 1, \pm 2, \cdots$.

For a hyponormal bilateral weighted shift $T$, we may assume all weight sequence of $T$ to have no zeros in the following discussion.

Now, the main results of this paper are as follows.

**Theorem 1** Let $T$ be a hyponormal bilateral weighted shift with $Te_n = \beta_n e_{n+1}$ for all $n = 0, \pm 1, \pm 2, \cdots$. If $0 < |\beta_n| < |\beta_{n+1}|$ for all $n = 0, \pm 1, \pm 2, \cdots$, then $T$ is near subnormal if and only if the sequence $\{\beta_n (\frac{|\beta_{n+1}|^2 - |\beta_n|^2}{|\beta_n|^2 - |\beta_{n-1}|^2})^{\frac{1}{2}}\}_{n=-\infty}^{+\infty}$ is bounded.

**Proof** Since $Te_n = \beta_n e_{n+1}$ for all $n = 0, \pm 1, \pm 2, \cdots$, it is easy to see that

$$Q_T = T^*T - TT^* = diag\{\cdots, |\beta_{-1}|^2 - |\beta_{-2}|^2, |\beta_0|^2 - |\beta_{-1}|^2, |\beta_1|^2 - |\beta_0|^2, \cdots\}. \tag{1}$$

Therefore,

$$Q_T^{\frac{1}{2}} = diag\{\cdots, (|\beta_{-1}|^2 - |\beta_{-2}|^2)^{\frac{1}{2}}, (|\beta_0|^2 - |\beta_{-1}|^2)^{\frac{1}{2}}, (|\beta_1|^2 - |\beta_0|^2)^{\frac{1}{2}}, \cdots\}, \tag{2}$$

and the M-P generalized inverse of $Q_T^{\frac{1}{2}}$ is

$$Q_T^{+\frac{1}{2}} = diag\{\cdots, \frac{1}{(|\beta_{-1}|^2 - |\beta_{-2}|^2)^{\frac{1}{2}}}, \frac{1}{(|\beta_0|^2 - |\beta_{-1}|^2)^{\frac{1}{2}}}, \frac{1}{(|\beta_1|^2 - |\beta_0|^2)^{\frac{1}{2}}}, \cdots\}. \tag{3}$$

By (2) and (3), through some calculations, it follows that $Q_T^{\frac{1}{2}} T Q_T^{+\frac{1}{2}}$ is also a bilateral weighted shift with weight sequence $\{\beta_n (\frac{|\beta_{n+1}|^2 - |\beta_n|^2}{|\beta_n|^2 - |\beta_{n-1}|^2})^{\frac{1}{2}}\}_{n=-\infty}^{+\infty}$. From (1), it is evident that $N(Q_T) = \{0\}$. By using the theorem of Section 1, we know that Theorem 1 holds.

Consider a hyponormal bilateral weighted shift $T$ with $Te_n = \beta_n e_{n+1}$ for all $n = 0, \pm 1, \pm 2, \cdots$. If $0 < |\beta_n| = |\beta_{n+1}|$ for all $n = 0, \pm 1, \pm 2, \cdots$, then $T$ is a normal operator.

For this reason, in the following discussion, we may exclude the above simple case and the case $0 < |\beta_n| < |\beta_{n+1}|$ for all $n = 0, \pm 1, \pm 2, \cdots$.

**Theorem 2** Let $T$ be a hyponormal bilateral weighted shift with $Te_n = \beta_n e_{n+1}$ for all $n$

=0, ±1,±2,⋯. If there exists an integer $j_0$ such that $|\beta_n| < |\beta_{n+1}|$ for all $n \leq j_0 - 1$, then $T$ is near subnormal if and only if there is an integer $k$ such that

( I ) $|\beta_n| < |\beta_{n+1}|$ for all $n \leq k - 1$, $|\beta_n| = |\beta_k|$ for all $n \geq k + 1$,

and ( II ) the sequence $\{\beta_n(\frac{|\beta_{n+1}|^2 - |\beta_n|^2}{|\beta_n|^2 - |\beta_{n-1}|^2})^{\frac{1}{2}}\}_{n=-\infty}^{n=k-1}$ is bounded.

**Proof** The proof of necessity. Suppose there exists an integer $j_0$ such that $|\beta_n| < |\beta_{n+1}|$ for all $n \leq j_0 - 1$. Now, we can choose $k$ to be the smallest integer such that $|\beta_k| = |\beta_{k+1}|$. Then $|\beta_n| < |\beta_{n+1}|$ for all $n \leq k - 1$, and $|\beta_k| = |\beta_{k+1}|$.

In this case, it is easy to see that $e_{k+1} \in N(Q_T)$. Since $T$ is near subnormal, so $N(Q_T) \in LatT$ and $Q_T^{\frac{1}{2}}TQ_T^{-\frac{1}{2}} \in B(H)$ by the theorem of Section1. It follows that $Te_{k+1} = \beta_{k+1}e_{k+2} \in N(Q_T)$. Hence $e_{k+2} \in N(Q_T)$ for $\beta_{k+1} \neq 0$. Similarly, we have $\{e_n\}_{n=k+1}^{+\infty} \subset N(Q_T)$. This fact implies that $|\beta_n| = |\beta_k|$ for all $n \geq k + 1$. Under present circumstances, the weight sequence of $Q_T^{\frac{1}{2}}TQ_T^{-\frac{1}{2}}$ is

$$\{\cdots, \beta_{k-2}(\frac{|\beta_{k-1}|^2 - |\beta_{k-2}|^2}{|\beta_{k-2}|^2 - |\beta_{k-3}|^2})^{\frac{1}{2}}, \beta_{k-1}(\frac{|\beta_k|^2 - |\beta_{k-1}|^2}{|\beta_{k-1}|^2 - |\beta_{k-2}|^2})^{\frac{1}{2}}, 0, 0, \cdots\}. \qquad (4)$$

From the theorem of Section1, it follows that the sequence $\{\beta_n(\frac{|\beta_{n+1}|^2 - |\beta_n|^2}{|\beta_n|^2 - |\beta_{n-1}|^2})^{\frac{1}{2}}\}_{n=-\infty}^{n=k-1}$ is bounded. Therefore, the necessity of Theorem 2 has been proved.

The proof of sufficiency. Suppose there is an integer $k$ such that $|\beta_n| < |\beta_{n+1}|$ for all $n \leq k - 1$, $|\beta_n| = |\beta_k|$ for all $n \geq k + 1$, and the sequence $\{\beta_n(\frac{|\beta_{n+1}|^2 - |\beta_n|^2}{|\beta_n|^2 - |\beta_{n-1}|^2})^{\frac{1}{2}}\}_{n=-\infty}^{n=k-1}$ is bounded. Then $Q_T = diag\{\cdots, |\beta_{k-1}|^2 - |\beta_{k-2}|^2, |\beta_k|^2 - |\beta_{k-1}|^2, 0, 0, \cdots\}$.

Hence $N(Q_T) = \overline{span\{e_{k+1}, e_{k+2}, \cdots\}}$. It is obvious that $N(Q_T) \in LatT$. Meanwhile, from (4) we know that $Q_T^{\frac{1}{2}} T Q_T^{+\frac{1}{2}} \in B(H)$. By the theorem of Section 1, $T$ is a near subnormal operator.

We have completed the proof of Theorem 2.

**Theorem 3** Let $T$ be a hyponormal bilateral weighted shift with $Te_n = \beta_n e_{n+1}$ for all $n$ $=0, \pm 1, \pm 2, \cdots$. If there exists an integer $j_0$ such that $|\beta_n| = |\beta_{j_0}|$ for all $n \leq j_0 - 1$, then $T$ is near subnormal if and only if $T$ is normal.

**Proof** The sufficiency is obvious. We turn to prove the necessity. If there exists an integer $j_0$ such that $|\beta_n| = |\beta_{j_0}|$ for all $n \leq j_0 - 1$, then

$$Q_T = T^*T - TT^* = diag\{\cdots, 0, 0, \cdots, 0, |\beta_{j_0+1}|^2 - |\beta_{j_0}|^2, |\beta_{j_0+2}|^2 - |\beta_{j_0+1}|^2, \cdots\}.$$

Therefore $e_{j_0} \in N(Q_T)$. Since $T$ is near subnormal, so $N(Q_T) \in LatT$. It follows that $Te_{j_0} = \beta_{j_0} e_{j_0+1} \in N(Q_T)$. Hence $e_{j_0+1} \in N(Q_T)$ for $\beta_{j_0} \neq 0$. Similarly, we have $\{e_n\}_{n=j_0}^{+\infty} \subset N(Q_T)$. This implies that $|\beta_n| = |\beta_{j_0}|$ for all $n \geq j_0$. Thus $T$ is normal.

Theorem 3 has been proved.

**Theorem 4** Let $T$ be a hyponormal bilateral weighted shift with $Te_n = \beta_n e_{n+1}$ for all $n$ $=0, \pm 1, \pm 2, \cdots$. If there exists an integer $j_0$ such that $|\beta_{j_0-1}| < |\beta_{j_0}| = |\beta_{j_0+1}| < |\beta_{j_0+2}|$, then $T$ is not near subnormal.

**Proof** Suppose there exists an integer $j_0$ such that $|\beta_{j_0-1}| < |\beta_{j_0}| = |\beta_{j_0+1}| < |\beta_{j_0+2}|$. Then $Q_T = T^*T - TT^* = diag\{\cdots, |\beta_{j_0}|^2 - |\beta_{j_0-1}|^2, 0, |\beta_{j_0+2}|^2 - |\beta_{j_0+1}|^2, \cdots\}$. Thus $e_{j_0+1} \in N(Q_T)$.

Since $Q_T(Te_{j_0+1}) = Q_T(\beta_{j_0+1} e_{j_0+2}) = \beta_{j_0+1}(|\beta_{j_0+2}|^2 - |\beta_{j_0+1}|^2)e_{j_0+2} \neq 0$, so $N(Q_T) \notin LatT$. By the theorem of Section 1, $T$ is not near subnormal.

We have completed the proof of Theorem 4.

## 3  An Application

Now, by Theorems 1 and 2, we can construct many bilateral weighted shifts which are near subnormal but not subnormal. Specially, by Theorems 3 and 4, it is easy to provide a great many bilateral weighted shifts which are hyponormal but not near subnormal. All of these operators are further answers to the Hilbert space problem 160.

Let $T$ be a hyponormal bilateral weighted shift with $Te_n = \beta_n e_{n+1}$ for all $n = 0, \pm 1, \pm 2, \cdots$. Three examples are as follows at the end of this paper.

**Example 1** Let $\beta_n = \dfrac{1}{n}$ for all $n \leq -1,$ and $\beta_n = 2$ for all $n \geq 0$. Then it is easy to show that the sequence $\{\beta_n(\dfrac{|\beta_{n+1}|^2 - |\beta_n|^2}{|\beta_n|^2 - |\beta_{n-1}|^2})^{\frac{1}{2}}\}_{n=-\infty}^{n=0}$ converges to 0 when $n \to -\infty$. Hence the sequence $\{\beta_n(\dfrac{|\beta_{n+1}|^2 - |\beta_n|^2}{|\beta_n|^2 - |\beta_{n-1}|^2})^{\frac{1}{2}}\}_{n=-\infty}^{n=0}$ is bounded. It follows from Theorem 2 that $T$ is near subnormal.

**Example 2** Choose $\beta_n = \dfrac{1}{n-1}$ for all $n \leq -1,$ $\beta_0 = \dfrac{2}{3}$ and $\beta_n = 2 - \dfrac{1}{n}$ for all $n \geq 1$. Then it is not difficult to prove that the sequence $\{\beta_n(\dfrac{|\beta_{n+1}|^2 - |\beta_n|^2}{|\beta_n|^2 - |\beta_{n-1}|^2})^{\frac{1}{2}}\}_{n=2}^{n=+\infty}$ converges to 2 when $n \to +\infty$. From Example 1 , the sequence $\{\beta_n(\dfrac{|\beta_{n+1}|^2 - |\beta_n|^2}{|\beta_n|^2 - |\beta_{n-1}|^2})^{\frac{1}{2}}\}_{n=-\infty}^{+\infty}$ is bounded. By Theorem 1, $T$ is near subnormal.

Moreover, we can show that the operators in Examples 1 and 2 are not subnormal by using the criteria in [1] or [6] for subnormal operators .We omit the details.

**Example 3** Choose two positive numbers $\lambda, \mu$ such that $\lambda < \mu$. Suppose $|\beta_n| = \lambda$ for all $n \leq 0,$ and $|\beta_n| = \mu$ for all $n \geq 1$. Then by Theorem 3, $T$ is hyponormal but not near subnormal .